\documentclass[a4paper,leqno,10pt,final]{amsproc}
\usepackage{amssymb,latexsym,amsxtra,amscd,amsfonts,amsthm,amsmath,verbatim}
\numberwithin{equation}{section}

\usepackage{epsfig}
\usepackage{amsfonts}
\usepackage{amsmath}

\usepackage{euscript}
\usepackage{amscd}

\def\and{\operatorname{and}}

\def\ker{\operatorname{ker}}
\def\im{\operatorname{im}}
\def\dim{\operatorname{dim}}
\def\depth{\operatorname{depth}}

\newcommand{\mm}{\mathfrak m}
\newcommand{\m}{\mathfrak m}

\newcommand{\ZZ}{\mathbb Z}
\newcommand{\NN}{\mathbb N}

\newcommand{\R}{\mathcal R}

\newcommand{\C}{\mathcal C}
\newcommand{\wrt}{ with respect to }

\newcommand{\bl}{\begin{lemma}}
\newcommand{\el}{\end{lemma}}
\newcommand{\bt}{\begin{theorem}}
\newcommand{\et}{\end{theorem}}
\newcommand{\ben}{\begin{enumerate}}
\newcommand{\een}{\end{enumerate}}
\newcommand{\bpf}{\begin{proof}}
\newcommand{\beqn}{\begin{eqnarray*}}
\newcommand{\eeqn}{\end{eqnarray*}}
\newcommand{\bd}{\begin{definition}}
\newcommand{\ed}{\end{definition}}
\newcommand{\bp}{\begin{proposition}}
\newcommand{\ep}{\end{proposition}}
\newcommand{\bc}{\begin{corollary}}
\newcommand{\ec}{\end{corollary}}
\newcommand{\bex}{\begin{example}}
\newcommand{\eex}{\end{example}}
\newcommand{\lra}{\longrightarrow}
\newcommand{\CM}{ Cohen-Macaulay }

\theoremstyle{plain}
\newtheorem{theorem}[equation]{Theorem}
\newtheorem{corollary}[equation]{Corollary}
\newtheorem{proposition}[equation]{Proposition}
\newtheorem{lemma}[equation]{Lemma}

\theoremstyle{definition}

\newtheorem{example}[equation]{Example}
\newtheorem{definition}[equation]{Definition}

\theoremstyle{remark}

\begin{document}

\noindent
{\em Advances in Algebra and Geometry}\\
{ \em Proceedings of the Hyderabad Conference, 2001}\\
{\em Hindustan Publishing Company}\\

%%% amsproc has separate commands for each of the following,
%%%whereas under article.cls
%%% (which was used by hba.sty), many of these went into \thanks
%%% also note that \thanks is not inside \author{} as it was in article.cls
\title{Local cohomology modules of bigraded Rees algebras}
\date{ } % leave blank

\author{A. V. Jayanthan \and J. K. Verma}
\thanks{\\ The first author is supported by the National Board for Higher
Mathematics, India 
}
\address{Department of Mathematics, Indian Institute of Technology
Bombay, Powai, Mumbai, India - 400076}
\email{jayan@math.iitb.ac.in}
\email{jkv@math.iitb.ac.in}
\maketitle

\centerline{\em Dedicated to Prof. R. C. Cowsik on the occasion of his
sixtieth birthday}

\begin{abstract}
Formulas are obtained in terms of complete reductions for the bigraded
components of local cohomology
modules of bigraded Rees algebras of 0-dimensional ideals in
2-dimensional Cohen-Macaulay local rings. As a consequence,
cohomological expressions for the coefficients of the Bhattacharya
polynomial of such ideals are obtained.
\end{abstract}

\section{Introduction}

Let $(R,\m)$ be a $d$-dimensional local ring. Let $I$ and $J$ be
$\m$-primary ideals of $R$. Let $\lambda$ denote length.
The function $B(r,s) = \lambda(R/I^rJ^s)$
is given by a polynomial $P(r,s)$ for large values of $r$ and $s$ \cite{b}.
  The
function $B(r,s)$ is called the {\em Bhattacharya function} of $I$ and $J$
and the polynomial $P(r,s)$ is called the {\em Bhattacharya polynomial} of
$I$ and $J$. The polynomial $P(r,s)$ has total degree $d$
\cite{b}. and it can be written as
$$
P(r,s) = \displaystyle{\sum_{0\leq i,j \leq d}e_{ij}{r \choose i}{s
\choose j}}.
$$
The integers $e_{ij}$ for which $i+j = d$ are called the {\em mixed
multiplicities} of $I$ and $J$.

Associated to $I$ and $J$ is the {\em bigraded Rees algebra}, $\R =
\oplus_{r,s\geq0}I^rJ^st_1^rt_2^s$. Here $t_1$ and $t_2$ are
indeterminates. Put $\R_{++} = It_1Jt_2\R$. Let $H^i_{\R_{++}}(\R)$
denote the local cohomology module of $\R$ with support in $\R_{++}$.

In this paper we study relationships among the Bhattacharya
function $B(r,s)$, the Bhattacharya polynomial $P(r,s)$ and local
cohomology modules $H^i_{\R_{++}}(\R)$ and obtain cohomological 
expressions for all the
coefficients of $P(r,s)$ when $R$ is a two-dimensional \CM local ring.

A relationship between Hilbert coefficients and local cohomology was first
observed by Grothendieck and Serre. Their formula expresses the
difference between Hilbert polynomial and the Hilbert function of a
graded algebra in terms of the lengths of graded components of local
cohomology modules with support in the maximal homogeneous ideal
\cite [Theorem 4.4.3]{bh}.

Sally, in \cite{sa1}, studied the Hilbert function of the maximal
ideal $\m$ of a Cohen-Macaulay local ring and its relation with the
local cohomology modules of the Rees algebra $R(\m) =
\oplus_{n\geq 0}\m^nt^n$ with support in $\m tR(\m)$.
Before we state her results, we recall the
concept of reductions of ideals introduced by Northcott and Rees in
\cite{nr}. An ideal $J$ of a Noetherian local ring $(R, \m)$ is said
to be a {\em reduction} of an ideal $I$ if $J \subseteq I$ and $JI^n =
I^{n+1}$ for $n \gg 0$. We say $J$ is a {\em minimal reduction} of $I$ if
$J$ is minimal with respect to inclusion among all reductions of $I$. If
$R/\m$ is infinite, then all minimal reductions of $I$ are generated
by same number of elements, $l(I)$, called the {\em analytic spread} of $I$.
The analytic spread of $I$ is equal to the dimension of the {\em fiber cone}
$F(I) = \oplus_{n\geq0}I^n/\m I^n$. We have ht $(I) \leq l(I)
\leq \min\{\mu(I), \dim(R)\}$, \cite{r1}. Here ht($I$) denotes the height
of $I$ and $\mu(I) = \dim_{R/\m}I/\m I$ denotes the minimum number of
generators for $I.$

Let ${\bf x} = x_1, \ldots, x_d \in \m$ be a minimal reduction
of $\m$. Let $({\bf x})^{[k]}$ denote the ideal $(x_1^k, \ldots,
x_d^k)$. For a graded module $M$ we will use the symbol $M_n$ or $[M]_n$
for  the $nth$ graded component of $M.$ Sally showed \cite{sa1} that
$$
\left[H^d_{R(\m)_+}(R(\m))\right]_n \simeq
\displaystyle\lim_{\stackrel{\longrightarrow}{k}}\m^{dk+n}/({\bf
x})^{[k]}\m^{(d-1)k+n}.
$$
Sally showed \cite{sa1} that the modules appearing on the right hand side
of the above equation are isomorphic for large $k$. The maps in the directed
system of modules appearing in the above direct limit are given by
multiplication
by $x_1x_2\ldots x_d.$  The above results were
generalized for $\m$-primary ideals, in dimension 2, by Sally
\cite{sa2}, and for {\em Hilbert filtration} of ideals in a $d$-dimensional
Cohen-Macaulay local ring by Blancafort in her thesis, \cite{bl}. By Hilbert
filtration of ideals $\{I_n\}$ we mean a descending sequence of
ideals of $(R,\m)$ such that $I_1I_n\subseteq I_{n+1}$ for all $n$,
$I_1I_n=I_{n+1}$ for large $n$, and $I_1$ is $\m$-primary.

It is natural to ask for analogues of these results for the bigraded
Rees algebra $\R$. We find the analogues in dimension $2.$ We have studied
the case of arbitrary dimension in \cite{jv} where completely different
methods are employed. In dimension two we can approach the results in a
self contained manner from the first principles. To state these analogues
we recall the concept of complete reduction introduced by Rees in
\cite{r2}. Let $(R,\m)$ be a Noetherian local ring of dimension $d$ and
$I$, $J$ be $\m$-primary ideals of $R$. Let $x_1, \ldots, x_d \in I$ and
$y_1, \ldots, y_d \in J$ and let $z_i = x_iy_i$ for $i = 1, \ldots, d$.
Then the system of elements $(x_1, \ldots, x_d, y_1, \ldots, y_d)$ is said
to be a complete reduction of $(I,J)$ if $(z_1, \ldots, z_d)$ is a
reduction of $IJ$. D. Rees showed that the complete reductions exist
if the local ring has infinite residue field.

Let $(R,\m)$ be a 2-dimensional local ring and $I, J$ be
$\m$-primary ideals of $R$. Let $(x,y;z,w)$ be a complete reduction of
$(I,J)$. We will show in section 4 that
$$
[H^2_{\R_{++}}(\R)]_{(r,s)} =
\displaystyle\lim_{\stackrel{\longrightarrow}{k}}
\frac{I^{2k+r}J^{2k+s}}{(xz, yw)^{[k]}I^{k+r}J^{k+s}}.
$$
We will also show that the modules appearing on the right hand side
of the above equation are isomorphic for $k \gg 0$.
Sally, \cite{sa2} proved,
$$
H^2_{B_+}(B)_n = \displaystyle\lim_{\stackrel{\longrightarrow}{k}}
I^{2k+n}/({\bf x})^{[k]}I^{k+n}
$$
where $B = R[It]$ and $B_+$ is the positively graded ideal of $B$.
She also showed that \cite[Proposition 2.4]{sa2} for all $n \geq 0,$
$$
\lambda(H^2_{B_+}(B)_n) = P_I(n) - \lambda(R/\widetilde{I^n}).
$$
Here $\tilde{I^n}$ denotes the {\em Ratliff-Rush  closure} of $I^n$ (see
section 2)  and  $P_I(n)$ denotes the Hilbert-Samuel polynomial
corresponding to the function $\lambda(R/I^n)$. Following \cite{sa2},
write $P_I(n)$ as
$$
P_I(n) = e(I){n \choose 2} + a_1(I) n + a_2(I).
$$
In Section 3 we will provide a simple and short proof, in dimension two, of
a
theorem
of Rees concerning mixed multiplicities and complete reductions \cite{r2}.
The purpose of providing
this proof is to make this paper self-contained.
In section 4 we will show that in dimension two, for $r, s \geq 0$,
$$
P(r,s) = \lambda\left([H^2_{\R_{++}}(\R)]_{(r,s)}\right) +
\lambda(R/\widetilde{I^rJ^s}).
$$

Blancafort obtained an expression for $H^1_{R(I)_+}(R(I))_n$ where $I$
is an $\m$-primary ideal of a Cohen-Macaulay local ring of dimension
at least 2. She showed that for $n \geq 0$
$$
H^1_{R(I)_+}(R(I))_n \cong \widetilde{I^n}/I^n.
$$
We prove, in section 4, that for $r, s \geq 0$,
$$
H^1_{\R_{++}}(\R)_{(r,s)} \cong \widetilde{I^rJ^s}/I^rJ^s.
$$

Sally obtained the following cohomological expressions for
the coefficients of the
Hilbert-Samuel polynomial $P_{\m}(n)$:

\begin{enumerate}
\item $\lambda\left(H^2_{R(\m)_+}(R(\m))_0\right) =  a_2(\m)$.
\item $\lambda\left(H^2_{R(\m)_+}(R(\m))_1\right) =  a_1(\m)+a_2(\m)-1$.
\item $\lambda\left(H^2_{R(\m)_+}(R(\m))_{-1}\right) = e(\m)-a_1(\m)+a_2(\m)$.
\end{enumerate}
The above formulas have also been deduced in \cite{jov}.
We obtain similar cohomological expressions for the the
Bhattacharya coefficients as a consequence of the expressions obtained
for $H^1_{\R_{++}}(\R)$ and $H^2_{\R_{++}}(\R)$. We  end the
article with an example to show that our results are no longer true if we
remove the Cohen-Macaulay hypothesis on the ring.
\vskip 2mm
\noindent
{\bf Acknowledgement:} The authors would like to thank the referee for
a very careful reading and for suggesting  several improvements and
additions.

\section{Ratliff-Rush closure of products of ideals}

Since the formulas for the local cohomology modules $H^i_{\R_{++}}(\R)$
involve Ratliff-Rush closure of products of ideals, we develop their basic
properties in this section.
Let $R$ be a Noetherian ring and $I$ be an ideal of $R$. The stable
value of the sequence $\{I^{n+1} : I^n\}$ is called the Ratliff-Rush
closure of $I$, denoted by $\tilde{I}$. An ideal $I$ is said to be
Ratliff-Rush if $\tilde{I} = I$.
The following proposition summarizes some basic
properties of Ratliff-Rush closure found in \cite{rr}.

\begin{proposition}\label{rr}
Let $I$ be an ideal containing a regular element in a Noetherian ring
$R$. Then
\begin{enumerate}
\item $I \subseteq \tilde{I}$ and $\widetilde{(\tilde{I})} =
\tilde{I}$.
\item $(\tilde{I})^n = I^n$ for $n \gg 0$. Hence if $I$ is
$\mm$-primary, the Hilbert polynomial of $I$ and $\tilde{I}$ are same.
\item $\widetilde{(I^n)} = I^n$ for $n \gg 0$.
\item If $(x_1, \ldots, x_g)$ is a reduction of $I$, then
$\tilde{I} = \cup_{n \geq 0} I^{n+1} : (x_1^n, \ldots, x_g^n)$.
\end{enumerate}
\end{proposition}

%In our description of local cohomology modules of bigraded Rees algebras
%we will use the notion of complete reductions of ideals introduced by
%D. Rees in \cite{r2}. We now recall this important concept.
%\vskip 2mm
%\noindent
%Let $(R,\m)$ be a Noetherian local ring of dimension $d$ and $I$, $J$
%be $\m$-primary ideals of $R$. Let $x_1, \ldots, x_d \in I$ and $y_1,
%\ldots, y_d \in J$ and let $z_i = x_iy_i$ for $i = 1, \ldots, d$. Then
%the system of elements $(x_1, \ldots, x_d, y_1, \ldots, y_d)$ is said
%to be a complete reduction of $(I,J)$ if $(z_1, \ldots, z_d)$ is a
%reduction of $IJ$. D. Rees showed that the complete reductions exist
%if the local ring has infinite residue field.
%In \cite{jv}, it is
%proved that the Ratliff-Rush closure of product of two ideals
%can be computed by taking colons with certain powers of a complete
%reduction.  For the sake of completeness we will include the proof of
%the next result.
In \cite{jv}, the authors studied some of the properties of Ratliff-Rush
closure of products of ideals. We recall some of the results proved there.
For the sake of completeness, we include the proof of the next result.
\vskip 2mm
\noindent

\begin{lemma}\label{rrprod}
Let $(R,\m)$ be a $d$-dimensional Noetherian local ring with infinite
residue field and $I, J$ ideals of $R$.  Then
\begin{enumerate}
\item $\widetilde{IJ} = \displaystyle{\bigcup_{r,s \geq 0}}
I^{r+1}J^{s+1} : I^rJ^s$.
\item $\widetilde{I^aJ^b} = \displaystyle{\bigcup_{k \geq
0}}I^{a+k}J^{b+k} : I^kJ^k$.
\item Let $I$ and $J$ are $\mm$-primary ideals with a
reduction $( {\bf z})=(z_1, \ldots, z_d)$ of $IJ$ obtained from a
complete
reduction of $I$ and $J$. Put $  {\bf z}^{[k]}=(z_1^k, \ldots, z_d^k)$ then
$$
\widetilde{I^aJ^b} = \displaystyle{\bigcup_{k \geq 0}}I^{a+k}J^{b+k}
:   {\bf z}^{[k]}.
$$
%\item If $IJ$ has a reduction generated by regular elements, then
%$\widetilde{I^rJ^s} = I^rJ^s$ for $r, s \gg 0$.
\end{enumerate}
\end{lemma}
\begin{proof}
(1). Let $x \in \widetilde{IJ}$, then $xI^nJ^n \subseteq I^{n+1}J^{n+1}$
for some $n$. Conversely if $xI^rJ^s
\subseteq I^{r+1}J^{s+1}$ for some $r,s \geq 0$ then for $n = max
\{r, s\}$, $xI^nJ^n \subseteq I^{n+1}J^{n+1}$ so that
$x \in (\widetilde{IJ})$. \\
(2). Let $z \in (\widetilde{I^aJ^b})$
then for some $r, s,$ by (1), we have $zI^{ar}J^{bs} \subseteq
I^{ar+a}J^{bs+b}$. Set
$k = max\{ar, bs\}$. Then $zI^kJ^k \subseteq I^{a+k}J^{k+b}$ and hence
$z \in I^{a+k}J^{b+k} : I^kJ^k$. Let $zI^kJ^k \subseteq
I^{a+k}J^{b+k}$ for some $k$. We may assume that $k
    = nab$ for $n \gg 0$. Therefore
    $z \in I^{nab+a}J^{nab+b} : I^{nab}J^{nab} \subseteq
(\widetilde{I^aJ^b})$. \\
(3). Suppose $y \in (\widetilde{I^aJ^b})$. Then for some $k$,
$yI^kJ^k \subseteq I^{a+k}J^{b+k}$, by (2). Since
$  {\bf z}^{[k]} \subseteq
I^kJ^k$, we have $y  {\bf z}^{[k]} \subseteq I^{a+k}J^{b+k}$.
Now, let $yz_i^k \in I^{a+k}J^{b+k}$ for $i = 1, \ldots, d$.
Write the reduction equation of $IJ$
\wrt $ {\bf z} $:
$   (IJ)^{m+n} =  {\bf z}^m(IJ)^n$ for all $m
\geq 0$ and $n$ large. Then, for $r$ large,
    $(IJ)^{r+dk} = ({\bf z })^{dk}I^rJ^r$.
Therefore,
$$
yI^{r+dk}J^{r+dk}  =  y( {\bf z})^{dk}I^rJ^r
     =  \sum_{\sum i_j = dk}y z_1^{i_1}\cdots z_d^{i_d}I^rJ^r
    \subseteq  I^{a+dk}J^{b+dk}I^rJ^r
$$
Hence $y \in (\widetilde{I^aJ^b})$, by (2). \\

\end{proof}

\section{Complete reductions and mixed multiplicities in dimension
$2$}
In the next section, while deriving a formula for the second local
cohomlogy module of the bigraded Rees algebra $\R$ of two $\m$-primary
ideals $I$ and $J$ in a two dimensional Cohen-Macaulay local ring $(R,\m)$,
    we will use a result of Rees linking mixed multiplicities
and complete reductions. Recall that a set of elements
$(x_1, x_2, \ldots, x_d)$ in a local ring  $(R,\m)$ of dimension $d$
is called a {\em joint reduction} of a set of $\m$-primary ideals
$(I_1, I_2, \ldots, I_d)$ if the ideal
    $x_1 I_2 I_3 \ldots I_d+\cdots+x_d I_1 I_2 \ldots I_{d-1}$ is a reduction
of the
product $I_1 I_2 \ldots I_d.$

Rees showed \cite{r3} that if
$(x_1,x_2, \ldots, x_d)$ is a joint reduction of
$(I_1, I_2, \ldots, I_d)$ where
$I_1=I_2=\cdots = I_i=I$ and $I_1=I_2=\cdots = I_j=J$, $i+j=d$, then
$e_{ij}=e(x_1, x_2, \ldots, x_d).$ The converse is true in quasi-unmixed
local rings and it was proved in dimension $2$ by Verma  \cite{sw1} and
in any dimension by Swanson in  \cite{sw2}. We will need the following
version for complete reductions.

\bt\label{rees}
Let $(R,\m)$ be a $2$-dimensional local ring and $I$ and $J$ be
$\m$-primary ideals. Let $(x,y;z,w)$ be a complete reduction of $(I,J).$
Then $$e(x,w)=e_{11}=e(y,z).$$
\et

\begin{proof}
First note that since $(xz,yw)$ is a reduction of $IJ$, there exists
an $n$ such that $(xz,yw)(IJ)^n = (IJ)^{n+1}.$ Hence

$$xI^nJ^{n+1}+wJ^nI^{n+1}=yI^nJ^{n+1}+zJ^nI^{n+1}=(IJ)^{n+1}.$$

Therefore $(x, w)$ and $(y, z)$ are both joint reductions of $(I, J).$
It is enough to prove that $e(y,z) = e_{11}.$ Consider the $R$-module
homomorphism
$$ \phi: R/I^n \oplus R/J^n \longrightarrow (y^n, z^n)/(y^nJ^n+ z^nI^n) $$
given by $\phi(\bar{r}, \bar{s})=rz^n+sy^n+(y^nJ^n+z^nI^n).$
Since $\phi$ is surjective, $$\lambda(R/y^nJ^n+z^nI^n)) \leq
\lambda(R/(y^n, z^n))+\lambda(R/I^n)+\lambda(R/J^n).$$
Divide by $n^2/2$ and take the limit as $n \rightarrow \infty,$ to get
$$
\lim_{n\rightarrow \infty}
\lambda\left(\frac{(R/(y^nJ^n+z^nI^n))}{n^2/2}\right) \leq
\lim_{n\rightarrow \infty}\frac{\lambda(R/(y^n,z^n))}{n^2/2} +e(I)+e(J).
$$

By Lech's Lemma \cite[Theorem 14.12]{m2}, the first term on the right
hand side of the above inequality is $2e(y,z)$ and  the left hand side is
bounded below by
  $$
\lim_{n \rightarrow \infty}\frac{\lambda\left(R/(yJ+zI)^n\right)}{n^2/2}=
e(yJ+zI)
.$$
  As $yJ+zI$ is a reduction of $IJ$, $e(IJ)=e(yJ+zI).$
Hence $e(IJ) \leq 2e(y,z)+e(I)+e(J).$ Since $e(IJ)=e(I)+2e_{11}+e(J),$
by the definition of the Bhattacharya polynomial, it follows that
$e_{11} \leq e(y,z).$ By repeated application of \cite[Section 5]{l},  we get
$$e(x,y)+e(x,w)+e(z,y)+e(z,w) = e(xz,yw) = e(IJ) = 2e_{11}+e(I)+e(J).$$
Hence $2e_{11}=e(x,w)+e(z,y).$ As $e_{11} \leq e(y,z)$ and
$e_{11} \leq e(x,w),$ we get  $e_{11}=e(x,w)=e(z,y).$
\end{proof}

\section{Local cohomology of bigraded Rees algebras}

For the rest of the article, let $(R,\m)$ denote a 2-dimensional
Cohen-Macaulay local ring with infinite residue field, unless stated
otherwise.
%For $x_1,
%\ldots, x_g \in R$, we denote by $({\bf x})^{[k]}$ the ideal $(x_1^k,
%\ldots, x_g^k)$.
Let $I$ and $J$ be  $\m $-primary ideals of $R$ and $(x,y;z,w)$
be a complete reduction of $I$ and $J$.

Let $\R = R[It_1, Jt_2]=\oplus_{r,s \in \NN }I^rJ^st_1^rt_2^s$ denote 
the bigraded
Rees algebra of $R$ with respect to $I$ and $J$. Let $\R(k, k)$ denote
the ring $\R$ with the $(r,s)$th graded piece $\R(k, k)_{(r,s)}
= \R_{(r+k, s+k)}$. We explain our method for computing the formulas
for the local cohomology modules of the Rees algebra $\R$.
    Consider the complex:

$$
{F^k}^. : 0 \longrightarrow \mathcal{R}
\buildrel\alpha_k\over\longrightarrow \mathcal{R}(k,k)^2
\buildrel\beta_k\over\longrightarrow \mathcal{R}(2k,2k) \longrightarrow
0,
$$
where the maps are defined as,
$$
\alpha_k(1) = ((xzt_1t_2)^k, (ywt_1t_2)^k)~~ \and~~ \beta_k(u, v) =
(ywt_1t_2)^ku - (xzt_1t_2)^kv.
$$
The twists are given so that $\alpha_k$ and $\beta_k$ are degree zero
maps.  We have the following commutative diagram:

$$
\CD
0@>   >> \R @> \alpha_{k}  >> \R(k,k)^2 @> \beta_{k} >>\R(2k,2k) @>  >> 0 \\
@.     @Vf_kVV             @Vg_kVV                @Vh_kVV   @.  \\
0@>   >> \R @> \alpha_{k+1}  >> \R(k+1,k+1)^2 @> 
\beta_{k+1} >>\R(2k+2,2k+2)@> >> 0
\endCD
$$

Here the map $f_k$ is the identity map, $g_k(1,0)=(xzt_1t_2,0),$
$g_k(0,1)=(0,ywt_1t_2)$ and $h_k(1)=xzyw(t_1t_2)^2.$ Thus the cohomology
modules $H^i(F_k^.)$ of the complex ${F_k}^.$ form a directed system for
each $i.$

Let $\C = (xz, yw)$ and $\C {\bf t}=(xzt_1t_2, ywt_1t_2)$.
Since both $\R_{++}$ and $\C{\bf t}$ have same radical and
$H^i_{\C{\bf t}}(\R) =
\displaystyle\lim_{\stackrel{\longrightarrow}{k}}
H^i({F^k}^{.})=H^i_{\R_{++}}(\R)$
for all
$i$ by  \cite[Theorem 5.2.9]{bs}.

    We begin by deriving an expression for the
bigraded components of the first local cohomology module of the
bigraded Rees algebra.

\bp\label{h1}
$[H^1_{\R_{++}}(\R)]_{(r,s)} \cong \widetilde{I^rJ^s}/I^rJ^s$ for all
$r, s \geq 0$.
\ep
\begin{proof}

By the above discussion we have:
$$
[H^1_{\R_{++}}(\R)]_{(r,s)} =
\displaystyle\lim_{\stackrel{\longrightarrow}{k}}
\frac{(\ker \beta_k)_{(r,s)}}{(\im \alpha_k)_{(r,s)}}.
$$

We show that for large $k$, $(\ker \beta_k)_{(r,s)} \cong
\widetilde{I^rJ^s}$ and $(\im \alpha_k)_{(r,s)} \cong I^rJ^s$.  Let
$(ut_1^{r+k}t_2^{s+k}, vt_1^{r+k}t_2^{s+k}) \in (\ker
\beta_k)_{(r,s)}$ for $k \gg 0$. Then $(yw)^ku -
(xz)^kv = 0$. Since $(xz, yw)$ is a
regular sequence in $R$, we have $u = p(xz)^k$ for some $p \in R$ and
$v = p(yw)^k$. Therefore
$$
(ut_1^{r+k}t_2^{s+k}, vt_1^{r+k}t_2^{s+k}) =
p((xz)^kt_1^{r+k}t_2^{s+k}, (yw)^kt_1^{r+k}t_2^{s+k}).
$$
Hence $(\ker \beta_k)_{(r,s)} \subseteq \{p((xz)^kt_1^{r+k}t_2^{s+k},
(yw)^kt_1^{r+k}t_2^{s+k}) \mid p \in R\}.$ The reverse inclusion is
clear.
%It remains to show that $(\ker \beta)_{(r,s)} = \widetilde{I^rJ^s}$.
%An element $p((xz)^kt_1^mt_2^n, (yw)^kt_1^mt_2^n) \in (\ker
%\beta)_{(r,s)}$ if and only if $p(xz)^k, p(yw)^k \in I^{k+r}J^{k+s}$.
Therefore
\begin{eqnarray*}
(\ker \beta_k)_{(r,s)} & = & \{p((xz)^kt_1^{r+k}t_2^{s+k},
(yw)^kt_1^{r+k}t_2^{s+k}) \mid p \in R\}.
\end{eqnarray*}
Consider the map $\delta: (\ker \beta_k)_{(r,s)} \longrightarrow
I^{r+k}J^{s+k} : (xz, yw)^{[k]}$ defined by
$$
\delta(p((xz)^kt_1^{r+k}t_2^{s+k}, (yw)^kt_1^{r+k}t_2^{s+k})) = p.
$$
It is clear that $\delta$ is an isomorphism. Hence, by Lemma
\ref{rrprod}(3), $\ker \beta_k(r,s) \cong \widetilde{I^rJ^s}$.

To see that $(\im \alpha_k)_{(r,s)} \cong I^rJ^s$, consider the map
$\phi : I^rJ^s \longrightarrow (\im \alpha_k)_{(r,s)}$ defined by
$\phi(p) = (p(xz)^kt_1^{r+k}t_2^{s+k},\; p(yw)^kt_1^{r+k}t_2^{s+k})$.
Since $xz$ and $yw$ are regular, $\phi$ is an isomorphism.

Therefore
$$
[H^1_{\R_{++}}(\R)]_{(r,s)} \cong \widetilde{I^rJ^s}/I^rJ^s.
$$

\end{proof}
\vskip 2mm
\noindent

Next we obtain an expression for the second local cohomology module of
bigraded Rees algebra in terms of complete reductions. In this result
we do not need the usual \CM hypothesis.

\bp\label{h2formula} Let $(R,\m)$ be a Noetherian local ring of dimension
$2$ and $I$, $J$ be $\m$-primary ideals of $R$. Let $(x, y; z, w)$ be a
complete reduction of $(I, J)$, where $x,y \in I$ and $z, w \in J$.  Set
$\R_{++} = (xzt_1t_2, ywt_1t_2)$. Then
$$
[H^2_{\R_{++}}(\R)]_{(r,s)} \cong
\displaystyle\lim_{\stackrel{\longrightarrow}{k}}
\frac{I^{2k+r}J^{2k+s}}{(xz, yw)^{[k]}I^{k+r}J^{k+s}}.
$$
\ep
\begin{proof} Consider the complex defined above
\begin{eqnarray*}
{F^{k}}^{.} : 0 \longrightarrow \mathcal{R}
\buildrel\alpha_k\over\longrightarrow \mathcal{R}(k,k)^2
\buildrel\beta_k\over\longrightarrow \mathcal{R}(2k,2k) \longrightarrow
0.
\end{eqnarray*}

Thus $H^2_{\R_{++}}(\R) =
\displaystyle\lim_{\stackrel{\longrightarrow}{k}}\mathcal{R}(2k,2k)/\im
\beta_k$. Note that the local cohomology modules have a natural
$\ZZ^2$-grading which is inherited from the $\NN^2$-grading of $\R$.
Therefore
$$
[H^2_{\R_{++}}(\R)]_{(r,s)} =
\displaystyle\lim_{\stackrel{\longrightarrow}{k}}
\frac{I^{2k+r}J^{2k+s}t_1^{2k+r}t_2^{2k+s}}{(\im \beta_k)_{(r,s)}}.
$$
Since $(\im \beta_k)_{(r,s)} =
(xzt_1t_2, ywt_1t_2)^{[k]}(It_1)^{k+r}(Jt_2)^{k+s},$
$$
[H^2_{\R_{++}}(\R)]_{(r,s)} \cong
\displaystyle\lim_{\stackrel{\longrightarrow}{k}}
\frac{I^{2k+r}J^{2k+s}}{(xz,yw)^{[k]}I^{k+r}J^{k+s}}.
$$
\end{proof}
For the directed system involved in the above direct limit, the maps
$$\frac{I^{2k+r}J^{2k+s}}{(xz,yw)^{[k]}I^{k+r}J^{k+s}}
\buildrel\mu\over\longrightarrow
\frac{I^{2k+r+2}J^{2k+s+2}}{(xz,yw)^{[k+1]}I^{k+r+1}J^{k+s+1}}$$
are the multiplication by $(xzyw)$. We show that the above map is
an isomorphism for $k \gg 0$ in the next two lemmas.
\begin{lemma}
The map $\mu$, defined as above, is surjective for large $k$.
\end{lemma}
\begin{proof}
To show that the maps are surjective for large $k$, we need to see
that for $k \gg 0$,
$$
I^{2k+r+2}J^{2k+s+2} \subseteq xzywI^{2k+r}J^{2k+s}
+  (xz, yw)^{[k+1]}I^{k+r+1}J^{k+s+1}
$$
Since $(xz,yw)$ is a reduction of $IJ$, for  $k \gg 0$,
\beqn
I^{2k+r+2}J^{2k+s+2} & = & (xz,
yw)^{k+1}I^{k+r+1}J^{k+s+1} \\
              & = & (xz, yw)^{[k+1]}I^{k+r+1}J^{k+s+1} \\
              & + & \sum_{i=1}^{k}((xz)^i
                    (yw)^{k+1-i})I^{k+r+1}J^{k+s+1}\\
              & \subseteq & xzywI^{2k+r}J^{2k+s}
	   + (xz, yw)^{[k+1]}I^{k+r+1}J^{k+s+1}.
\eeqn
Hence $\mu$ is surjective.
\end{proof}
%\vskip 2mm
%\noindent
%We need the following lemma to prove the injectivity of $\mu$.
%\bl\label{ker}
%Let $r, s \geq 0$ and $k \geq 2$ and let
%$$
%\phi : \frac{I^{2+r}J^{2+s}}{(xz, yw)I^{1+r}J^{1+s}}
%\longrightarrow \frac{I^{2k+r}J^{2k+s}}{(xz,
%yw)^{[k]}I^{k+r}J^{k+s}}
%$$
%be the multiplication by $(xzyw)^{k-1}$. Then
%$$
%\ker \phi \subseteq \frac{(xz, yw)\widetilde{(I^{r+1}J^{s+1})}}
%{(xz, yw)(I^{r+1}J^{s+1})}.
%$$
%\el
%\begin{proof} Let $d \in \ker \phi$. Then there exist $p ,q  \in
%I^{k+r}J^{k+s}$ such that $(xzyw)^{k-1}d = (xz)^kp  + (yw)^kq $.  Then
%$(xz)^{k-1}((yw)^{k-1}d - xz p) = (yw)^{k}q$ Since $(xz, yw)$ is a
%regular sequence in $R$, $q = (xz)^{k-1}u$ and $p = (yw)^{k-1}v$ for
%some $u, v \in R$. Substituting the values of $p$ and $q$ and
%canceling $(xz)^{k-1}(yw)^{k-1}$ from both sides, we get $d = xzu +
%ywv$ with $u \in I^{k+r+1}J^{k+s+1} : (xz)^k$ and $v \in
%I^{k+r+1}J^{k+s+1} : (yw)^k$. Therefore to complete the proof, it is
%enough to show that $u, v \in \widetilde{I^{r+1}J^{s+1}}$.
%
%\noindent
%{\it Claim :} $u \in I^{k+r+1}J^{k+s+1} : (yw)^k$ and $v \in
%I^{k+r+1}J^{k+s+1} : (xz)^k$. We have
%\beqn
%(xz)^{k-1}(yw)^{k-1}d & = & (xz)^kp + (yw)^kq \\
                         %& = & (xz)^kp + (yw)^k(xz)^{k-1}u \\
%\eeqn
%Therefore $(yw)^{k-1}d - xzp  =  (yw)^ku$, so that $u \in
%I^{k+r+1}J^{k+s+1} : (yw)^k$. Similarly one can show that $v \in
%I^{k+r+1}J^{k+s+1} : (xz)^k$.  It follows that $u, v \in
%I^{k+r+1}J^{k+s+1} :(xz, yw)^{[k]} \subseteq
%\widetilde{I^{r+1}J^{s+1}}$, by Lemma \ref{rrprod}.
%\end{proof}
\vskip 2mm
\noindent
\bl\label{h2poly}
Fix $r,s$. With the notation as above, the multiplication map $\mu =
\mu_{xzyw}$
$$
\frac {I^{2k+r}J^{2k+s}}{(xz,
yw)^{[k]}I^{k+r}J^{k+s}}\buildrel\mu\over\longrightarrow
\frac{I^{2k+r+2}J^{2k+s+2}}{(xz,yw)^{[k+1]}I^{k+r+1}J^{k+s+1}}
$$
is an isomorphism for $k \gg 0$.
\el

\begin{proof} Let $d \in I^{2k+r}J^{2k+s}$ and suppose that $\bar{d}
\in \ker \mu$. Then
\begin{eqnarray}\label{eqn1}
xzywd = (xz)^{k+1}p + (yw)^{k+1}q
\end{eqnarray}
with $p, q \in I^{k+r+1}J^{k+s+1}.$ Then $xz(ywd - (xz)^kp) =
(yw)^{k+1}q$. Since $(xz, yw)$ is a regular sequence, $q \in (xz)$.
Let $q = uxz$ for some $u \in R$. Similarly, $p = vyw$ for some $v \in
R$. Therefore $u \in I^{k+r+1}J^{k+s+1} : xz$ and $v \in
I^{k+r+1}J^{k+s+1} : yw$ and hence $u \in I^{2k+r}J^{2k+s} : (xz)^k$
and $v \in I^{2k+r}J^{2k+s} : (yw)^k$ for $k \gg 0$. We  show that $u,
v \in \widetilde{I^{k+r}J^{k+s}}$.  Substituting the value of $q$ in
(\ref{eqn1}) and cancelling $xz$ on both sides we get, $ywd = (xz)^kp +
(yw)^{k+1}u$. Therefore $(yw)^{k+1}u = ywd - (xz)^kp \in
I^{2k+r+2}J^{2k+r+2}$. Thus $u \in I^{2k+r+1}J^{2k+r+1} : (yw)^{k+1}$.
Hence, for $k \gg 0$, $u \in I^{2k+r+1}J^{2k+r+1} : ((xz)^{k+1},
(yw)^{k+1}) = \widetilde{I^{k+r}J^{k+s}}$, by Lemma \ref{rrprod}(3).
Since $IJ$ has a reduction
generated by regular elements, we can apply Lemma 3.4 of \cite{jv} to see
that for $k \gg 0$, $(\widetilde{I^{k+r}J^{k+s}}) = (I^{k+r}J^{k+s})$.
Thus the map is injective for large $k$. We have already shown that
the map is surjective for large $k$. Therefore $\mu$ is an isomorphism
for $k \gg 0$.
\end{proof}
\vskip 2mm
\noindent
Lemma \ref{h2poly} shows that to compute
$\lambda([H^2_{\R_{++}}(\R)]_{(r,s)})$
it enough to compute the length of
$$\frac{I^{2k+r}J^{2k+s}}{(xz,yw)^{[k]}I^{k+r}J^{k+s}} ~~~
\mbox{ for } k \gg 0.$$
We aim to compute $\lambda\left([H^2_{\R_{++}}(\R)]_{(r,s)}\right)$
for which we need the following technical result.
\bl
With notations as before, for $k \gg 0$,
$$
\lambda \left(\frac{(xz, yw)^{[k]}}{(xz, yw)^{[k]}I^{k+r}J^{k+s}}\right)
= 2 \lambda\left(\frac{R}{I^{k+r}J^{k+s}}\right) -
\lambda(R/\widetilde{I^rJ^s}).
$$
\el
\begin{proof} Consider the exact sequence
$$
0 \longrightarrow K \longrightarrow
\left(\frac{R}{I^{k+r}J^{k+s}}\right)^2
\buildrel\alpha\over\longrightarrow \frac{(xz, yw)^{[k]}}{(xz,
yw)^{[k]}I^{k+r}J^{k+s}} \longrightarrow 0
$$
where $\alpha(\overline{g},
\overline{h}) = \overline{(xz)^kg + (yw)^kh} $ and  $K = \ker \alpha$.
First we compute $K$. Let  $(\overline{g}, \overline{h}) \in K$. Then
there exist $p, q \in I^{k+r}J^{k+s}$ such that
$$ (xz)^kg+ (yw)^kh = (xz)^kp + (yw)^kq \hspace*{0.2in} \cdots ~(2).$$
i.e  $ (xz)^k(p -g) = (yw)^k(h -q)$. Thus $p -g \in ((yw)^k)$, say $p
- g = (yw)^ku$ for some $u \in R$. Substituting in (2) and cancelling
$(yw)^k$ we get, $g = p
- (yw)^ku$ and $h = q + (xz)^ku$. Therefore $(\overline{g},
\overline{h}) = u(-\overline{(yw)^k}, \overline{(xz)^k})$. Therefore
$K = \ker \alpha = (-\overline{(yw)^k},
\overline{(xz)^k})R$. Define $\phi : R \longrightarrow K$
by $\phi(u) = u(-\overline{(yw)^k}, \overline{(xz)^k})$. Clearly
$\phi$ is surjective. Also
$$
\ker \phi = \{u \in R : u(yw)^k, u(xz)^k \in
I^{k+r}J^{k+s}\} = I^{k+r}J^{k+s} : (xz, yw)^{[k]}.
$$
For $k \gg 0, ~ I^{k+r}J^{k+s} : (xz, yw)^{[k]} =
(\widetilde{I^rJ^s})$. Thus $K \cong R/(\widetilde{I^rJ^s})$ for $k
\gg 0$. Hence for large $k$, $$\lambda\left(\frac{(xz, yw)^{[k]}}{(xz,
yw)^{[k]}I^{k+r}J^{k+s}}\right) = 2\lambda(R/I^{k+r}J^{k+s}) -
\lambda(R/(\widetilde{I^rJ^s})).$$
\end{proof}

\begin{theorem}\label{h2}
Let $r,s \geq 0$. With notations as before
$$
\lambda\left((H^2_{\R_{++}}(R))_{(r,s)}\right) = P(r,s) -
\lambda(R/\widetilde{I^rJ^s}).
$$
\end{theorem}
\begin{proof}
By Proposition \ref{h2formula} and Lemma \ref{h2poly}, for $k \gg 0$
\beqn
\lambda\left((H^2_{\R_{++}}(\R))_{(r,s)}\right)
     & = &
\lambda\left(\frac{I^{2k+r}J^{2k+s}}{(xz,yw)^{[k]}I^{k+r}J^{k+s}}\right) \\
     & = &
\lambda\left(\frac{(xz, yw)^{[k]}}{(xz,yw)^{[k]}I^{k+r}J^{k+s}}\right)
+ \lambda\left(\frac{R}{(xz,yw)^{[k]}}\right) \\
     &  &
    - ~~\lambda\left(\frac{R}{I^{2k+r}J^{2k+s}}\right) \\
     & = &
2 \lambda\left(\frac{R}{I^{k+r}J^{k+s}}\right) -
\lambda\left(\frac{R}{(\widetilde{I^rJ^s})}\right) -
\lambda\left(\frac{R}{I^{2k+r}J^{2k+s}}\right)  \\
     &  &
    + ~~ \lambda\left(\frac{R}{(xz,yw)^{[k]}}\right).
\eeqn

Hence for $k \gg 0$, we have
\beqn
\lambda\left((H^2_{\R_{++}}(\R))_{(r,s)}\right) & = & 2P(k+r, k+s) -
P(2k+r,2k+s) \\
     & + &
\lambda\left(\frac{R}{(xz,yw)^{[k]}}\right)
-\lambda\left(\frac{R}{\widetilde{I^rJ^s}}\right) ~~\cdots(*)
\eeqn
%\end{proof}
Since $(x,y)$, $(x,w)$, $(z,y)$ and $(z,w)$ are regular sequences in
$R$, we get
\beqn
\lambda\left(\frac{R}{(xz,yw)^{[k]}}\right)
     & = & k^2\lambda(R/(xz,yw)) \\
     & = & k^2[\lambda(R/(x,y)) + \lambda(R/(x,w)) + \lambda(R/(z,y)) +
     \lambda(R/(z,w))]
\eeqn
so that, by Theorem \ref{rees}
\beqn
\lambda\left((H^2_{\R_{++}}(\R))_{(r,s)}\right)
     & = & 2P(k+r, k+s) - P(2k+r,2k+s)  \\
     & + & k^2(e(I) + 2e_{11} + e(J)) -
     \lambda\left(\frac{R}{(\widetilde{I^rJ^s})}\right).
\eeqn

Substituting for $P(r,s)$ and expanding we get
\begin{eqnarray}\label{h2formula1}
\lambda\left((H^2_{\R_{++}}(R))_{(r,s)}\right) = P(r,s) -
\lambda(R/\widetilde{I^rJ^s}).
\end{eqnarray}
\end{proof}
The formulas obtained for $H^1_{\R_{++}}(\R)$ and $H^2_{\R_{++}}(\R)$
enable us to express the difference $P(r,s) - B(r,s)$ in terms of the
Euler characteristic of local cohomology of $\R$. The formula proved
below has been generalized to arbitrary dimension in \cite{jv}. The
proof in arbitrary dimension is quite different from the proof given
here. The proof of the general result uses analysis of cohomology of
modified Koszul complex.

Set $h^i(r,s) =
\lambda([H_{\R_{++}}^i(\R)]_{(r,s)})$.

\bt\label{2d-diffform}
$P(r,s) - B(r,s) = \sum_{i=0}^2
(-1)^ih^i(r,s)$ for all $r,s \geq 0$.
\et
\begin{proof}
Since $\R_{++}$ contains a regular element,  $H^0_{\R_{++}}(\R) = 0$
so that $h^0(r,s) = 0$. By Theorem \ref{h1} we have
$h^1(r,s) = \lambda(\widetilde{I^rJ^s}/I^rJ^s)$ and by
Theorem \ref{h2} we have $h^2(r,s) = P(r,s) - \lambda
(R/(\widetilde{I^rJ^s}))$.  Therefore
\beqn
\sum_{i=0}^2 (-1)^ih^i(r,s)
     & = &
P(r,s) - \lambda (R/(\widetilde{I^rJ^s}))
-\lambda(\widetilde{I^rJ^s}/I^rJ^s) \\
     & = & P(r,s) - B(r,s).
\eeqn
\end{proof}

One does not know {\em apriori} the stage at which the Bhattacharya
function equals the Bhattacharya polynomial. Due to this,
it is difficult to calculate the coefficients of the Bhattacharya
polynomials. Therefore it is desirable to have effective methods for
computing these coefficients since they contain information about
the bigraded Rees algebra, for example, its depth \cite{jv}. The formulas
we have obtained for the bigraded components of the local cohomology
modules of bigraded Rees algebra make it possible to find formulas
for the Bhattacharya coefficients.

\begin{corollary}\label{coefforms} With notations as before
\begin{enumerate}
\item $e_{00} = h^2(0,0)$.
\item $e_{10} = h^2(1,0) - h^2(0,0) + \lambda(R/\tilde{I}).$
\item $e_{01} = h^2(0,1) - h^2(0,0) + \lambda(R/\tilde{J}).$
\item $e_{20} = h^2(2,0) - 2h^2(1,0)+h^2(0,0) - 2\lambda(R/\tilde{I})
+ \lambda(R/\widetilde{I^2}).$
\item $e_{02} = h^2(0,2) - 2h^2(0,1)+h^2(0,0) - 2\lambda(R/\tilde{J})
+ \lambda(R/\widetilde{J^2}).$
\end{enumerate}
\end{corollary}
\begin{proof} (1) Putting $(r,s) = (0,0)$ in (\ref{h2formula1}) we get,
$e_{00} = h^2(0,0).$
\vskip 2mm
\noindent
(2) Put $(r,s) = (1,0)$ in (\ref{h2formula1}) to get $e_{10} + e_{00}
= h^2(1,0) + \lambda(R/\tilde{I}).$ This yields (2).
\vskip 2mm
\noindent
(3) Now put $(r,s) = (2,0)$ in (\ref{h2formula1}) to get $e_{20} +
2e_{10} + e_{00} = h^2(2,0) + \lambda(R/\widetilde{I^2}).$
Substituting the values of $e_{10}$ and $e_{00}$, we get
$$
e_{20} = h^2(2,0) - 2h^2(1,0)+h^2(0,0) - 2\lambda(R/\tilde{I})
+ \lambda(R/\widetilde{I^2}).
$$
Similarly one can get the required expressions for $e_{01}$ and $e_{02}$.
\end{proof}
%We denote $\lambda(H^i_{\R_{++}}(\R)_{(r,s)})$ by
%$h^i_{\R_{++}}(r,s)$.

We conclude the paper with an example to show that our results need
the Cohen-Macaulay hypothesis.

\bex Let $k$ be a field and $T=k[[X,Y,Z]]$ be a power series
ring over $k.$
The ring $R=T/(X^2,XY)T$ is a two-dimensional local ring of depth one. Let
$x$(resp. $y$ and $z$) be images of $X$(resp. $Y$ and $Z$) in $R.$  Then
$x^2=xy=0.$
Put $S=k[X,Y,Z]$ and  $L=(X^2,XY)S.$
Then $L=(X)S\cap (X^2, Y)S.$ Put $I=(X)S$ and $J=(X^2, Y)S.$
Let $\m$ denote the unique maximal ideal of $R.$ Then
the associated graded ring $G:= G(\m)\simeq S/L.$
Indeed, let $R^{\ast}$ denote the $\m$-adic completion of $R.$
Then $R^{\ast}\simeq T^{\ast}/(X^2,XY)=S/L$ by \cite[Theorems 54 and 55]{m}.
Since $G(\m) \simeq G(\m^{\ast}),$ and $R^{\ast}$ is homogeneous, it
follows that  $G \simeq S/L.$
To find the Hilbert series of $G$,
  consider the  exact sequence:
$$0\lra G \stackrel{\alpha}{\lra} S/I \oplus S/J
   \stackrel{\beta}{\lra} S/(I+J) \lra 0.$$
Here $\alpha(r')=(r', r')$ and $\beta((a',b'))=(a-b)'.$
Hence,

$$
H(G,t):=\sum_{n=0}^{\infty}\lambda\left(\m^n/\m^{n+1}t^n\right)t^n
= H(S/I, t) + H(S/J, t)-H(S/(I+J), t).
$$
Since
$$
H(S/I, t)=1/(1-t)^2, \;
H(S/J, t)=(1+t)/(1-t)\; \mbox{ and }
H(S/(I+J), t)=1/(1-t),
$$
we see that
$$H(G,t)=\frac{1+t-t^2}{(1-t)^2}.$$
Write the Hilbert polynomial $P(n)$ corresponding to the Hilbert
function $\lambda(R/\m^n)$ in the following way:
$$
P(n) = e(\m){n+1 \choose 2} - e_1(\m) n + e_2(\m).
$$
Then, by \cite[Proposition 4.1.9]{bh},
$e(\m)=1, e_1(\m)=-1$ and $e_2(\m)=-1.$
Thus  $P(n)$  is given by
$$
P(n) =  {n+1 \choose 2} + n -1.
$$

Therefore
\begin{eqnarray*}
\lambda(R/\m^r\m^s) & = & {r+s+1 \choose 2} + r+s -1 \\
& = & {r \choose 2} + rs + {s \choose 2} + 2(r+s) - 1.
\end{eqnarray*}
This shows that $e_{00} = -1$ which is not equal to the length of any
module. Therefore Theorem 4.7 (and hence Corollary 4.10) does not hold
if the assumption of Cohen-Macaulayness on the ring is removed.

We now show that for $I=J=\m,$ \; $H^2_{\R_{++}}(\R)=0.$
Since $\dim G(\m)=2,$ any  reduction minimally  generated by $2$
elements  is a minimal reduction \cite{nr}. Since $x^2=0,$ the
ideal $(y,z)$  is a reduction of $\m.$ Hence it is a minimal reduction
of $\m.$
%Since $x^2=0,$ the ideal $(y,z)$ is a minimal reduction of
%$\m.$
By Proposition \ref{h2formula},
$$
[H^2_{\R_{++}}(\R)]_{(r,s)} =
\displaystyle\lim_{\stackrel{\longrightarrow}{k}}
\frac{\m^{4k+r+s}}{(y^2, z^2)^{[k]}\m^{2k+r+s}}.
$$
For $n \geq 2 $
$$
\m^n=(y^n, y^{n-1}z, y^{n-2}z^2, \ldots, y^2z^{n-2},(x,y)z^{n-1},z^n).$$
%Since $(y,z)$ is a minimal reduction of $\m$,
It is easy to see, by using
the above formula for $\m^n$, that $(y^{2k},z^{2k})\m^{2k}=\m^{4k}.$
Therefore $(y^{2k},z^{2k})$ is a minimal reduction of $\m^{2k}.$
Hence by the above formula for $H^2_{\R_{++}}(\R),$ we conclude that
$H^2_{\R_{++}}(\R)=0.$ Since $\depth G(\m)=1,$
$\widetilde{\m^n}=\m^n$ for all $n\geq 0.$ This shows that
none of the formulas hold true in Corollary 4.10.

\eex

\end{document}